%% file: MasterThesis.tex
\documentclass{amsart}
	\title[A Destabilizing Wall of Skyscraper Sheaves on Ruled Surfaces]{Explicit Description of A Certain Destabilizing Wall of Skyscraper Sheaves on Ruled Surfaces}
	\author[Takayuki Uchiba]{Takayuki Uchiba}
	\address{Department of Pure and Applied Mathematics, Faculty of Science and Engineering, Waseda University, 3-4-1 Ohkubo Shinjyuku Tokyo 169-8555, Japan}
	\email{prime-xp.1990@moegi.waseda.jp}
\usepackage[varg]{txfonts}
\usepackage[dvips]{graphicx}
\usepackage{amscd}
	\DeclareMathOperator{\Hom}{Hom}
	\DeclareMathOperator{\Stab}{Stab}
	\DeclareMathOperator{\Coh}{Coh}
	\DeclareMathOperator{\rank}{rank}
	\DeclareMathOperator{\per}{per}
	\DeclareMathOperator{\re}{Re}
	\DeclareMathOperator{\NS}{NS}
	\DeclareMathOperator{\Amp}{Amp}
	\DeclareMathOperator{\ch}{ch}
	\DeclareMathOperator{\pr}{pr}
	\DeclareMathOperator{\td}{td}
	\DeclareMathOperator{\im}{Im}
	\DeclareMathOperator{\id}{id}

\begin{document}
\begin{abstract}
We give a explicit description of gluing stability conditions on ruled surfaces by introducing gluing perversity. Moreover, we describe a destabilizing wall of skyscraper sheaves on ruled surfaces by deformation of stability conditions glued from $\widetilde{GL^{+}}(2,\mathbb{R})$-translates of the standard stability condition on the base curve.
\end{abstract}

\maketitle

\section{Introduction}
Bridgeland introduced the notion of a \emph{stability condition} on a triangulated category in \cite{Bri1}. A \emph{stability space} which is a set of stability conditions on a fixed triangulated category has a natural topology if one assumes \emph{locally ﬁniteness} for stability conditions. Especially, each connected component of the stability space is a complex manifold (\cite{Bri1} Theorem 1.2). In this paper we describe a \emph{destabilizing wall} of skyscraper sheaves on ruled surfaces in the stability space. A fundamental example of locally ﬁnite stability condition is \emph{geometric stability conditions} (\cite{Bri2}\S6, \cite{Ohk} Deﬁnition 3.5). However, the skyscraper sheaves are stable of the same phase with respect to geometric stability conditions (\cite{Ohk} Proposition 3.6). Hence we need to ask if there is a stability condition with respect to which skyscraper sheaves are strictly semistable of the same phase.\\

Collins and Polishchuck \cite{CP} introduced \emph{gluing stability conditions} on a triangulated category that has a semi-orthogonal decomposition. A derived category on a ruled surface has a semi-orthogonal decomposition that consists of its subcategories which are equivalent to the derived category on the base curve (\cite{Orl}). Hence, one can hope to construct stability conditions glued from stability conditions on the base curve. In section 3, we introduce \emph{gluing perversity} (Deﬁnition 3.6), which is the key notion to the following theorem:
\newtheorem*{theorem11}{Theorem 1.1}
\begin{theorem11}[\protect Theorem 3.9]
On ruled surfaces, a stability condition $\sigma$ glued from $\widetilde{GL^{+}}(2,\mathbb{R})$-translates of the standard stability condition on the base curve is a locally ﬁnite stability condition if and only if the gluing perversity of $\sigma$ is at least one. 
\end{theorem11}
In this paper, we mean a stability condition glued from $\widetilde{GL^{+}}(2,\mathbb{R})$-translates of the \emph{standard stability condition} on the base curve simply by a gluing staiblity condition. One can see from Theorem 1.1 that the existence of gluing stability conditions does not depend on genus of ruled surfaces. This means that the gluing stability conditions constitute a class of fundamental stability conditions on ruled surfaces. Furthermore, we describe the following lemma on the stability of skyscraper sheaves in the description of gluing perversity.
\newtheorem*{lem12}{Lemma 1.2}
\begin{lem12}[\protect Lemma 3.10]
 Suppose that $\sigma$ is a gluing stability condition on a ruled surface.\\ 
{\rm(1)} If the gluing perversity of $\sigma$ is equal to 1, the skyscraper sheaves are strictly semistable of the same phase for any point of the ruled surface in $\sigma$.\\
{\rm(2)} If the gluing perversity is larger than 1, the skyscraper sheaves are not stable in for any point of the ruled surface in $\sigma$.
\end{lem12}
In section 4, we describe a destabilizing wall of skyscraper sheaves on ruled surfaces. Lemma 1.2 already suggests that the set of gluing stability conditions with gluing perversity 1 is a destabilizing wall in the stability space. By deformation theory of stability conditions (see \cite{Bri1} \S7.), we can prove the following lemma.
\newtheorem*{lem1.3}{Lemma 1.3}
\begin{lem1.3}[\protect From Lemma 4.2]
 Let $S$ be a ruled surface. Suppose that $\sigma_{gl}=(Z_{gl},P_{gl})$ is a gluing stability condition with the gluing perversity 1 on $S$. Then there is an $\epsilon_{0}>0$ such that if $0<\epsilon<\epsilon_{0}$ and $W:\mathcal{N}(S)\rightarrow\mathbb{C}$ is a group homomorphism satisfying\\
\centerline{the phase of $O_{f}(-C_{0})$ is greater than the phase of $O_{f}$, and}\\
\centerline{$|W(E)-Z(E)|<\sin(\pi\epsilon)|Z(F)|$}\\
for any $E\in D^{b}(S)$ semistable in $\sigma_{gl}$, then there is a unique locally finite Bridgeland stability condition $\tau=(W,\mathcal{Q})$ on $S$ with $d(\mathcal{P}_{gl},\mathcal{Q})<\epsilon$ satisfying that $O_{x}$ are stable of the same phase in $\tau$ for any $x\in S$.
\end{lem1.3}
From the above results, we can describe a certain destabilizing wall of skyscraper sheaves by simple calculation.
\newtheorem*{thm1.4}{Theorem 1.4}
\begin{thm1.4}[\protect From Theorem 4.4]
Let $S_{geom}$ be the set of geometric stability conditions on $S$ and $S_{gl,p}$ be the set of gluing stability conditions with gluing perversity p. Suppose that $A=(\begin{pmatrix}a&&\frac{1}{2}a\deg\mathcal{E}\\0&&a\end{pmatrix}^{-1},f)\in\widetilde{GL^{+}}(2,\mathbb{R})$ with $a<0$. Then $\partial\overline{S_{geom}}\cap S_{gl,1}$ is the set of $\widetilde{GL^{+}}(2,\mathbb{R})$-translates of a stability condition glued from $\sigma_{st}.A$ and $\sigma_{st}$. 
\end{thm1.4}
\newtheorem*{ack}{Acknowledgement}
\begin{ack}
The author would like to thank his surpervisor Prof. Yasunari Nagai for providing a great deal of help and continuously warm encouragement. He would also like to thank Prof. Yukinobu Toda, Prof. Hajime Kaji, Prof. Ryo Ohkawa, Prof. Daizo Ishikawa and  Prof. Seung-Jo Jung for many valuable comments, discussion and pointing out my mistakes. 
\end{ack}

\section{Geometric stability conditions on ruled surfaces}
Bridgeland introduced the notion of a stability condition on a triangulated category in \cite{Bri1}. 
\newtheorem*{def2.1}{Definition 2.1} 
\begin{def2.1}[\protect \cite{Bri1} Definition 5.1]
Let $\mathcal{D}$ be a trianguleted category and $K(\mathcal{D})$ Grothendieck group of $\mathcal{D}$. A \emph{Bridgeland stability condition} on $\sigma=(Z,\mathcal{P})$ on $\mathcal{D}$ consists of a linear map $Z:K(\mathcal{D})\rightarrow\mathbb{C}$ called the \emph{central charge}, and full additive subcategories $\mathcal{P}(\phi)\subset\mathcal{D}$ for each $\phi\in\mathbb{R}$, satisfying the following axioms.\\
{\rm(1)} for all $0\neq E\in \mathcal{P}(\phi)$, if there exists some $m(E)>0$ such that $Z(E)=m(E)\exp(i\pi\phi)$,\\
{\rm(2)} $\mathcal{P}(\phi+1)=\mathcal{P}(\phi)[1]$, for all $\phi\in\mathbb{R}$,\\
{\rm(3)} if $\phi_{1}>\phi_{2}$ and $A_{j}\in \mathcal{P}(\phi_{j})$ $(j=1,2)$ then $\Hom(A_{1},A_{2})=0$,\\
{\rm(4)} for each nonzero object $E\in D$, there is a finite sequence of real number\\
\centerline{$\phi_{1}>\phi_{2}>\cdot\cdot\cdot>\phi_{n}$}\\
and a collection of triangles\\
\centerline{$E_{j}\rightarrow E_{j+1}\rightarrow A_{j+1}\rightarrow E_{j}[1]$}\\
with $E_{0}=0$, $E_{n}=E$, and $A_{j+1}\in \mathcal{P}(\phi_{j+1})$ for all $j=0,\cdot\cdot\cdot,n-1$.
\end{def2.1}
$\mathcal{P}$ is called the \emph{slicing} of $\mathcal{D}$. An object $E$ is defined to be \emph{semistable} of phase $\phi$ in $\sigma$ if $E\in \mathcal{P}(\phi)$. A semistable object $E\in \mathcal{P}(\phi)$ is \emph{stable} if it has no nontrivial subobject in $\mathcal{P}(\phi)$.
\newtheorem*{def2.2}{Definition 2.2}
\begin{def2.2}[\protect \cite{Bri1} Definition 5.7]
A slicing $\mathcal{P}$ of a triangulated category $\mathcal{D}$ is \emph{locally finite} if there exists a real number $\eta>0$ such that the quasi abelian category $\mathcal{P}((t-\eta,t+\eta))\subset\mathcal{D}$ is of finite length for all $t\in\mathbb{R}$. A Bridgeland stability condition $(Z,\mathcal{P})$ is \emph{locally finite} if the corresponding slicing $\mathcal{P}$ is.
\end{def2.2}
Since the decomposition of a nonzero object $E\in\mathcal{D}$ given by Definition 2.1 (4) is unique up to isomorphisms, we can define $\phi^{+}_{\sigma}(E)=\phi_{n}$, $\phi^{-}_{\sigma}(E)=\phi_{1}$ and $m_{\sigma}(E)=\Sigma_{j}|Z(A_{j})|$. There is a generalized metric on the space of locally finite stability conditions $\Stab\mathcal{D}$ on a triangulated category $\mathcal{D}$. The metric $d$ is defined by\\
\centerline{$d(\sigma,\tau)=\sup_{0\neq E\in D}\left\{|\phi_{\sigma}^{+}(E)-\phi_{\tau}^{+}(E)|,|\phi_{\sigma}^{-}(E)-\phi_{\tau}^{-}(E)|,|\log\frac{m_{\tau}(E)}{m_{\sigma}(E)}|\right\}$.}\\
Then $\phi^{\pm}$ and $m(E)$ are continuous functions on $\Stab\mathcal{D}$. It follows immediately from this that the subset of $\Stab\mathcal{D}$ consisting of those stability conditions in which a given object is semistable is a closed subset (\cite{Bri1} Proposition 8.1).\\

Let $S$ be a smooth projective surface over $\mathbb{C}$. A Bridgeland stability condition $\sigma=(Z,\mathcal{P})$ is \emph{numerical} if the central charge $Z:K(S)\rightarrow\mathbb{C}$ factors through the numerical Grothendieck group $\mathcal{N}(S)$. \emph{Mukai pairing} is a symmetric bilinear form $(-,-)_{S}$ on $\mathcal{N}(S)\simeq\mathbb{Z}\oplus\NS(S)\oplus\frac{1}{2}\mathbb{Z}$ defined by the following formula\\
\centerline{$((r_{1},D_{1},s_{1}),(r_{2},D_{2},s_{2}))_{S}=D_{1}.D_{2}-r_{1}s_{2}-r_{2}s_{1}.$}\\
The set of numerical locally finite stability conditions $\Stab_{\mathcal{N}}S$ is called \emph{stability space}. If $\sigma=(Z,\mathcal{P})\in\Stab_{\mathcal{N}}S$, we can write $Z(E)=(\pr_{1}(\sigma),\ch(E))_{S}$. 
\newtheorem*{prop2.3}{Proposition 2.3}
\begin{prop2.3}[\protect \cite{Bri1} Corollary 1.3]
For each connected component $\Stab^{\dagger}S$\\$\subset\Stab_{\mathcal{N}}S$, there are a subspace $V(\Stab^{\dagger}S)\subset\Hom(\mathcal{N}(S),\mathbb{C})$ and a local homeomorphism $\pr_{1}:\Stab^{\dagger}S\rightarrow V(\Stab^{\dagger}S)$ which maps a stability condition to its central charge $Z$. In particular $\Stab^{\dagger}S$ is a finite dimensional complex manifold.
\end{prop2.3}
A connected component $\Stab^{\dagger}S$ is \emph{full} if the subspace $V(\Stab^{\dagger}S)$ is equal to $\Hom(\mathcal{N}(S),\mathbb{C})$. A stability condition $\sigma\in\Stab_{\mathcal{N}}S$ is \emph{full} if it lies in a full component. On a derived category of coherent sheaves on a surface, one of fundamental examples of numerical locally finite stability conditions are \emph{divisorial stability conditions} (\cite{AB} \S2). We can construct a divisorial stability condition in the following way:
\newtheorem*{def2.4}{Definition 2.4}
\begin{def2.4}[\protect \cite{Bri1} Definition 2.1]
Let $\mathcal{A}$ be an abelian category and $K(\mathcal{A})$ Grothendieck group of $\mathcal{A}$. A \emph{stability function} on $\mathcal{A}$ is a group homomorphism $Z:K(\mathcal{A})\rightarrow\mathbb{C}$ such that for all $0\neq E\in\mathcal{A}$ the complex number $Z(E)$ lies in the strict upper half plane $H=\left\{r\exp(i\pi\phi)\mid r>0\text{ and }0<\phi\leq1\right\}$.
\end{def2.4}
Let $\mathcal{A}$ be a heart of a bounded t-structure of a triangulated category $\mathcal{D}$. $\mathcal{A}$ is an abelian subcategory of $\mathcal{D}$ and one has an identification of Grothendieck group $K(\mathcal{D})=K(\mathcal{A})$. To give a stability condition on $\mathcal{D}$ is  equivalent to giving a bounded t-structure $\mathcal{D}$ and a stability function on its heart $\mathcal{A}$ with the Harder Narasimhan property (\cite{Bri1} Proposition 5.3).In this paper, stability function is also called \emph{pre-stability condition}.\\

We denote $\Amp(S)$ ample cone of $S$ and $\NS(S)$ Neron Severi group of $S$. Let $\omega\in\Amp(S)$. One defines the slope $\mu_{\omega}$ of a torsion free sheaf $E\in\Coh S$ by\\
\centerline{$\mu_{\omega}(E)=\frac{c_{1}(E).\omega}{\rank(E)}$.}\\
For any $B,\omega\in\NS(S)\otimes\mathbb{R}$ with $\omega\in\Amp(S)$ there is a unique torsion pair $(\mathcal{T}_{B,\omega},\mathcal{F}_{B,\omega})$ on the category $\Coh S$ such that $\mathcal{T}_{B,\omega}$ consists of sheaves whose torsion free parts have $\mu_{\omega}$-semistable Harder Narasimhan factors of slope $\mu_{\omega}>B.\omega$ and $\mathcal{F}_{B,\omega}$ consists of torsion free sheaves on $S$ all of whose $\mu_{\omega}$-semistable Harder Narasimhan factors have slope $\mu_{\omega}\leq B.\omega$ (\cite{Bri2} Lemma 6.1).
\newtheorem*{def2.5}{Definition 2.5}
\begin{def2.5}[\protect \cite{AB} \S2 Our Charges, \cite{Ohk} Definition 3.3]
$\sigma_{B,\omega}=(Z_{B,\omega},\mathcal{A}_{B,\omega})$ is defined by the stability function \\
\centerline{$Z_{B,\omega}(E)=(\exp(B+i\omega),\ch(E))_{S}$}\\
and the heart of the bounded t-structure $\mathcal{A}_{B,\omega}$ which is obtained from $\Coh S$ by tilting with respect to  the torsion pair $(\mathcal{T}_{B,\omega}, \mathcal{F}_{B,\omega})$.
\end{def2.5}
If $\sigma_{B,\omega}$ is a stability condition, Arcara and Miles called it a \emph{divisorial stability condition} (\cite{AM}). For each pair $B,\omega\in\NS(S)\otimes\mathbb{Q}$ with $\omega\in\Amp(S)$, $\sigma_{B,\omega}$ is a numerical locally finite stability condition (\cite{Ohk} Proposition 3.4). $\Stab\mathcal{D}$ carries a right action of the group $\widetilde{GL^{+}}(2,\mathbb{R})$, the universal covering space of $GL^{+}(2,\mathbb{R})$ (\cite{Bri2} Lemma 8.2). Ohkawa called a $\widetilde{GL^{+}}(2,\mathbb{R})$-translate of a divisorial stability condition a \emph{geometric stability condition} (\cite{Ohk} Definition 3.5).
\newtheorem*{prop2.6}{Proposition 2.6}
\begin{prop2.6}[\protect \cite{Ohk} Proposition 3.6]
$\sigma\in\Stab_{\mathcal{N}}S$ is geometric if and only if\\
{\rm(1)}  for all $x\in S$, skyscraper sheaves $\mathcal{O}_{x}$ are stable of the same phase in $\sigma$,\\
{\rm(2)} there exist $M\in GL^{+}(2,\mathbb{R})$ and $B,\omega\in\NS(S)\otimes\mathbb{R}$ such that $\omega^{2}>0$ and $M^{-1}\pr_{1}(\sigma)=\exp(B+i\omega)$.
\end{prop2.6}
A \emph{ruled surface} is a smooth projective surface $S$, together with a surjective morphism $p:S\rightarrow C$ to a smooth projective curve of genus $g$, such that the fibre $S_{x}$ is isomorphic to $\mathbb{P}^{1}$ for any point $x\in C$, and such that $p$ admits a section $s:C\rightarrow S$ (\cite{Har} \S V.2). Furthermore, let $C_{0}$ be $s(C)$, $\mathcal{E}$ the direct image sheaf $p_{*}\mathcal{O}_{S}(C_{0})$ and $f$ a fibre of $p$. Then $S$ is isomorphic to the projective bundle $\mathbb{P}_{C}(\mathcal{E})$ of $\mathcal{E}$, and we can calculate the intersection numbers as $C_{0}^{2}=\deg\mathcal{E}$, $C_{0}.f=1$, $f^{2}=0$, and the canonical divisor $K_{S}=-2C_{0}+(2g-2+\deg\mathcal{E})f$. $\NS(S)$ is generated by $C_{0}$ and $f$, and hence $\dim_{\mathbb{R}}\Hom(\mathcal{N}(S),\mathbb{C})=8$.
\newtheorem*{prop2.7}{Proposition 2.7}
\begin{prop2.7}
For any $f$, $\mathcal{O}_{f}$ is stable of the same phase in a geometric stability condition.
\end{prop2.7}
\begin{proof}{(c.f. \cite{Bri2} Lemma 6.3)}
 For any $f$, $Z(\mathcal{O}_{f})$ always take in the same value in $\mathbb{C}$. It follows immediately that for any $f$ the phase of $\mathcal{O}_{f}$ is the same if $O_{f}$ is stable. First, we show that a subobject of torsion sheaf is also torsion sheaf in $\mathcal{A}_{B,\omega}$. Suppose $T$ is torsion sheaf. Recall that $T$ lies in the torsion subcategory $\mathcal{T}_{B,\omega}$ and hence in the abelian category $\mathcal{A}_{B,\omega}$. Suppose that\\
\centerline{$0\rightarrow A\rightarrow T \rightarrow B\rightarrow 0$}\\
is a short exact sequence in $\mathcal{A}_{B,\omega}$ with $A\in\mathcal{T}_{B,\omega}$. Taking cohomology gives an exact sequence in $\Coh S$\\
\centerline{$0\rightarrow\mathcal{H}^{-1}(B)\rightarrow\mathcal{H}^{0}(A)\rightarrow T \rightarrow\mathcal{H}^{0}(B)\rightarrow 0$.}\\
Since $\mathcal{H}^{-1}(B)\in\mathcal{T}_{B,\omega}$, $\mathcal{H}^{-1}(B)$ is torsion free sheaf. It follows that the $\mu_{\omega}$-semistable facotrs of  $\mathcal{H}^{-1}(B)$ and $\mathcal{H}^{0}(A)$ have the same slope. The contradicts the definition of the category $\mathcal{A}_{B,\omega}$ unless $\mathcal{H}^{-1}(B)=0$, in which case either $A$ and $B$ must be torsion sheaf.\\
\ Second, we show that subobjects of $\mathcal{O}_{f}$ are $\mathcal{O}_{f}(-p_{1}-\cdot\cdot\cdot-p_{n})$ with $p_{1},\cdot\cdot\cdot,p_{n}\in f$. Let $i:f\hookrightarrow S$ and $F$ a subobject of $\mathcal{O}_{f}$. Then $F$ is a torsion sheaf and hence $i^{*}F$ is a subsheaf of the structure sheaf of $f$, which is $\mathcal{O}_{f}(-p_{1}-\cdot\cdot\cdot-p_{n})$ with $p_{1},\cdot\cdot\cdot,p_{n}\in f$. It follows that $F\simeq Ri_{*}i^{*}F=\mathcal{O}_{f}(-p_{1}-\cdot\cdot\cdot-p_{n})$ with $p_{1},\cdot\cdot\cdot,p_{n}\in f$. Hence, $\mathcal{O}_{f}$ is stable by comparison of these phases.\\
\end{proof}

\section{Constructing gluing stability conditions on ruled surfaces}
This section is concerned with the construction and the existence of the gluing stability conditions on ruled surfaces, and the stability of skyscraper sheaves in gluing stability conditions.\\

Since $p$ is a flat morphism, $p^{*}$ is an exact functor, and hence $Lp^{*}$ can be simply denoted by $p^{*}$. Since $\mathcal{O}_{S}(-C_{0})$ is locally free sheaf, $\otimes ^{L}\mathcal{O}_{S}(-C_{0})$ is ordinary tensor product $\otimes \mathcal{O}_{S}(-C_{0})$. Orlov \cite{Orl} showed that a derived category of a ruled surface has \emph{Orlov's semi-orthogonal decomposition} $D^{b}(S)=\langle p^{*}D^{b}(C)\otimes \mathcal{O}_{S}(-C_{0}),p^{*}D^{b}(C)\rangle$. Recall that $p^{*}D^{b}(C)\otimes \mathcal{O}_{S}(-C_{0})$ and $p^{*}D^{b}(C)$ are equivalent to triangulated category $D^{b}(C)$. There exist the following canonical isomorphisms of Grothendieck groups (c.f. \cite{MMS} section 2),\\
\centerline{$F_{1}:K(C)\simeq K(p^{*}D^{b}(C)\otimes \mathcal{O}_{S}(-C_{0}))$,}\\
\centerline{$F_{2}:K(C)\simeq K(p^{*}D^{b}(C))$.}\\
Furthermore, we can describe the space of stability conditions on the both categories,\\
\centerline{$\Stab(p^{*}D^{b}(C)\otimes \mathcal{O}_{S}(-C_{0}))=\left\{(Z_{1},\mathcal{P}_{1})\middle| \begin{array}{l}\text{$(Z,\mathcal{P})\in\Stab C$, $Z_{1}=Z\circ F_{1}^{-1}$}\\ \text{for all $\phi\in\mathbb{R}$ $\mathcal{P}_{1}(\phi)=p^{*}\mathcal{P}(\phi)\otimes \mathcal{O}_{S}(-C_{0})$} \end{array}\right\},$}\\
\centerline{$\Stab(p^{*}D^{b}(C))=\left\{(Z_{2},\mathcal{P}_{2})\middle| \begin{array}{l}\text{$(Z,\mathcal{P})\in\Stab C$, $Z_{2}=Z\circ F_{2}^{-1}$}\\ \text{for all $\phi\in\mathbb{R}$ $\mathcal{P}_{2}(\phi)=p^{*}\mathcal{P}(\phi)$} \end{array}\right\}.$}\\
$\Stab C$ is completely determined in \cite{Bri1}, \cite{Mac} and \cite{Oka}. $\sigma_{st}=(Z_{st},\mathcal{P}_{st})$ with $Z_{st}(E)=-\deg E+i \rank E$ and $\mathcal{P}(0,1]=\Coh C$ is a stability condition on $\Stab C$. It is called  \emph{standard stability condition}. Especially, the following result is remarkable.
\newtheorem*{prop3.1}{Proposition 3.1}
\begin{prop3.1}[\protect \cite{Bri1} Theorem 9.1, \cite{Mac} Theorem 2.7]
If a smooth projective curve $C$ has positive genus, then the action of $\widetilde{GL^{+}}(2,\mathbb{R})$ on $\Stab C$ is free and transitive, so that\\
\centerline{$\Stab C\simeq\widetilde{GL^{+}}(2,\mathbb{R})$.}
\end{prop3.1}
Collins and Polishchuck \cite{CP} gave the definition of gluing stability conditions.
\newtheorem*{def3.2}{Definition 3.2}
\begin{def3.2}[\protect \cite{CP} \S2. Definition]
Suppose $\mathcal{D}$ is a triangulated category that have a semi-orthogonal decomposition $\langle \mathcal{D}_{1},\mathcal{D}_{2}\rangle$, $\lambda_{1}$ is the left adjoint functor of $\mathcal{D}_{1}\rightarrow \mathcal{D}$ and $\rho_{2}$ is the right adjoint functor of $\mathcal{D}_{2}\rightarrow \mathcal{D}$. $\sigma=(Z,\mathcal{A})$ is called \emph{gluing pre-stability condition} of $\sigma_{1}$ and $\sigma_{2}$ if $\sigma_{j}=(Z_{j},\mathcal{A}_{j})\in\Stab\mathcal{D}_{j}$ ($j=1,2$) satisfy the following conditions,\\
{\rm(1)} $Z=Z_{1}\circ\lambda_{1}+Z_{2}\circ\rho_{2}$,\\
{\rm(2)} $\mathcal{A}=\left\{F\in \mathcal{D} \mid \lambda_{1}(\mathcal{F})\in \mathcal{A}_{1}\text{ and }\rho_{2}(\mathcal{F})\in \mathcal{A}_{2}\right\}$,\\
{\rm(3)} $\Hom(\mathcal{A}_{1},\mathcal{A}_{2}[i])=0$ for any $i\leq0$ (We call this \emph{gluing property}.)
\end{def3.2}
It is called \emph{gluing stability condition} if it satisfies Harder-Narasimhan property. In the above definition, we set $\mathcal{D}=D^{b}(S)$, $\mathcal{D}_{1}=p^{*}D^{b}(C)\otimes\mathcal{O}_{S}(-C_{0})$ and $\mathcal{D}_{2}=p^{*}D^{b}(C)$. Then we get explicit formulas of $\lambda_{1}$ and $\rho_{2}$.
\newtheorem*{prop3.3}{Proposition 3.3}
\begin{prop3.3}
Let $F$ be an object of $D^{b}(S)$. We get\\
{\rm(1)} $\lambda_{1}(F)=p^{*}(Rp_{*}(F(-C_{0}+(2g-2+\deg\mathcal{E})f))\otimes\omega_{C}^{*}[1])\otimes \mathcal{O}_{S}(-C_{0})$,\\
{\rm(2)} $\rho_{2}(F)=p^{*}Rp_{*}F$.
\end{prop3.3}
\begin{proof}
Recall that $p^{*}$ and $\otimes \mathcal{O}_{S}(-C_{0})$ are fully faithful. $\lambda_{1}$ can be calculated by the following calculation,\\
　$\Hom(F,p^{*}G\otimes \mathcal{O}_{S}(-C_{0}))$\\
　$=\Hom(F(C_{0}),p^{*}G)$\\
　$=\Hom(F(C_{0}), p^{!}G\otimes \omega_{p}^{*}[-1])$\\
　$=\Hom(F(C_{0})\otimes \omega_{p}[1],p^{!}G)$\\
　$=\Hom(Rp_{*}(F(C_{0})\otimes \omega_{p}[1]),G)$\\
　$=\Hom(Rp_{*}(F(C_{0})\otimes \omega_{S} \otimes p^{*}\omega_{C}^{*}[1]),G)$\\
　$=\Hom(p^{*}(Rp_{*}(F(-C_{0}+(2g-2+\deg\mathcal{E})f) \otimes p^{*}\omega_{C}^{*}[1]))\otimes \mathcal{O}_{S}(-C_{0}),p^{*}G\otimes \mathcal{O}_{S}(-C_{0}))$\\
　$=\Hom(p^{*}(Rp_{*}(F(-C_{0}+(2g-2+\deg\mathcal{E})f)) \otimes \omega_{C}^{*}[1])\otimes \mathcal{O}_{S}(-C_{0}),p^{*}G\otimes \mathcal{O}_{S}(-C_{0}))$\\
We can get $\rho_{2}$ by similar calculation. 
\end{proof}
If one takes stability conditions on $\mathcal{D}_1$ and $\mathcal{D}_2$, the gluing of the stability conditions under the above definition is not a stability condition. Gluing procedure is compatible with the action of $\widetilde{GL^{+}}(2,\mathbb{R})$.
\newtheorem*{prop3.4}{Proposition 3.4}
\begin{prop3.4}
Suppose $A\in\widetilde{GL^{+}}(2,\mathbb{R})$ and $\sigma_{gl}$ is a gluing pre-stability condition of $\sigma_{1}$ and $\sigma_{2}$. Then $\sigma_{gl}.A$ is equal to the gluing of $\sigma_{1}.A$ and $\sigma_{2}.A$.
\end{prop3.4}
\begin{proof}
By Definition 3.2 (2), both gluing stability conditions have the same central charge. We show that both have the same heart of the bounded t-structure. Let $A=(M,f)\in\widetilde{GL^{+}}(2,\mathbb{R})$. Suppose that $\sigma_{gl}=(Z_{gl},\mathcal{P}_{gl})$ is a stability condition glued from $\sigma_{1}=(Z_{1},\mathcal{P}_{1})$ and $\sigma_{2}=(Z_{2},\mathcal{P}_{2})$. For any $\phi$, $\mathcal{P}_{1}(f^{-1}(\phi))\subset\mathcal{P}_{gl}(f^{-1}(\phi))$ and $\mathcal{P}_{2}(f^{-1}(\phi))\subset\mathcal{P}_{gl}(f^{-1}(\phi))$ by \cite{CP} Proposition 2.2 (3). Then $\mathcal{P}_{1}(f^{-1}(0),f^{-1}(1)]\subset\mathcal{P}_{gl}(f^{-1}(0),f^{-1}(1)]$ and $\mathcal{P}_{2}(f^{-1}(0),f^{-1}(1)]$\\$\subset\mathcal{P}_{gl}(f^{-1}(0),f^{-1}(1)]$. Furthermore, we get the inclusion\\
\centerline{$\langle\mathcal{P}_{1}(f^{-1}(0),f^{-1}(1)],\mathcal{P}_{2}(f^{-1}(0),f^{-1}(1)]\rangle\subset\mathcal{P}_{gl}(f^{-1}(0),f^{-1}(1)]$}\\
by extension closedness. Hence, both have the same a heart of a bounded t-structure.
\end{proof}
 From now on, let $\sigma_{1}$ and $\sigma_{2}$ $\widetilde{GL^{+}}(2,\mathbb{R})$-translates of a stability condition on $p^{*}D^{b}(C)\otimes\mathcal{O}(-C_{0})$ and $p^{*}D^{b}(C)$ induced from the standard stability condition $D^{b}(C)$ respectively. We can calculate a central charge of such a gluing pre-stability condition.
\newtheorem*{prop3.5}{Proposition 3.5}
\begin{prop3.5}
Let $M^{-1}=\begin{pmatrix}a&b\\c&d\end{pmatrix}\in GL^{+}(2,\mathbb{R})$. Suppose that $\sigma_{1}$ is a stability condition on $p^{*}D^{b}(C)\otimes\mathcal{O}_{S}(-C_{0})$ and $\sigma_{2}$ is a standard stability condition on $p^{*}D^{b}(C)$. Then a gluing stability conditions $\sigma_{gl}=(Z_{gl},\mathcal{P}_{gl})$ glued from $\sigma_{1}.M$ and $\sigma_{2}$ satisfies\\
\centerline{$\pr_{1}(Z_{gl})=((1-a)-i c, -C_{0}+[\{\frac{1}{2}\deg\mathcal{E}(a+1)-b\}+i \{\frac{1}{2}c\deg\mathcal{E}+(1-d)\}]f,-i)$.}
\end{prop3.5}
\begin{proof}
 By Definition 3.2 (2) and Proposition 3.3, all we need to calculate is $\ch Rp_{*}(F(-C_{0}+(2g-2+\deg\mathcal{E})f))\otimes\omega_{C}^{*}[1]$ and $\ch Rp_{*}(F)$.\\
　Now, we calculate $\ch Rp_{*}(F(-C_{0}+(2g-2+\deg\mathcal{E})f))\otimes\omega_{C}^{*}[1]$. By Grothendieck-Riemann-Roch formula,\\
　$\ch Rp_{*}(F(-C_{0}+(2g-2+\deg\mathcal{E})f))\otimes\omega_{C}^{*}[1]$\\
　$=-p_{*}(\ch F(-C_{0}+(2g-2+\deg\mathcal{E})f).\td S).\td C^{-1}.\ch\omega_{C}^{-1}$.\\
Suppose that $\ch F=(r,c_{1},\ch_{2})$, then\\
　$\ch F(-C_{0}+(2g-2+\deg\mathcal{E})f)$\\
　$=\ch F.\ch\mathcal{O}_{S}(-C_{0}+(2g-2+\deg\mathcal{E})f)$\\
　$=(r,c_{1}-rC_{0}+r(2g-2+\deg\mathcal{E})f,\ch_{2}-c_{1}.C_{0}+(2g-2+\deg\mathcal{E})c_{1}.f+\frac{1}{2}r(-4g+4-\deg\mathcal{E}))$.\\
　$\ch F(-C_{0}+(2g-2+\deg\mathcal{E})f).\td S$\\
　$=(r,c_{1}+\frac{1}{2}r(2g-2+\deg\mathcal{E})f,\ch_{2}+\frac{1}{2}(2g-2+\deg\mathcal{E})c_{1}.f)$.\\
　$p_{*}(\ch F(-C_{0}+(2g-2+\deg\mathcal{E})f).\td S).\td C^{-1}$\\
　$=(c_{1}.f,\ch_{2}+(2g-2+(\frac{1}{2}\deg\mathcal{E})c_{1}.f))$.\\
　$p_{*}(\ch F(-C_{0}+(2g-2+\deg\mathcal{E})f).\td S).\td C^{-1}.\ch\omega_{C}^{-1}$\\
　$=(c_{1}.f,\ch_{2}+(\frac{1}{2}\deg\mathcal{E})c_{1}.f)$.\\
Hence, $\ch Rp_{*} F(-C_{0}+(2g-2+\deg\mathcal{E})f)\otimes\omega_{C}^{*}[1]=(-c_{1}.f,-\ch_{2}-(\frac{1}{2}\deg\mathcal{E})c_{1}.f)$. We can get $\ch Rp_{*}(F)=(c_{1}.f+r,\ch_{2}+c_{1}.C_{0}-(\frac{1}{2}\deg\mathcal{E})c_{1}.f)$ similarly. Then we get\\
　$\re Z_{gl}(F)$\\
　$=[a\{\ch_{2}+(\frac{1}{2}\deg\mathcal{E})c_{1}.f\}+b(-c_{1}.f)]-\{\ch_{2}+c_{1}.C_{0}-(\frac{1}{2}\deg\mathcal{E})c_{1}.f\}$\\
　$=-c_{1}.C_{0}+\{\frac{1}{2}\deg\mathcal{E}(a+1)-b\}c_{1}.f+(a-1)\ch_{2}$\\
　$\im Z_{gl}(F)$\\
　$=[c\{\ch_{2}+(\frac{1}{2}\deg\mathcal{E})c_{1}.f\}+d(-c_{1}.f)]+(c_{1}.f+r)$\\
　$=r+\{(\frac{1}{2}c\deg\mathcal{E})c_{1}.f+(1-d)c_{1}.f\}+c\ch_{2}$
\end{proof}
Now, one cannot usually glue $\sigma_1$ and $\sigma_2$. For describing a necessary and sufficient condition of the existence of the gluing stability condition, we introduce \emph{gluing perversity}.
\newtheorem*{def3.6}{Definition 3.6}
\begin{def3.6}
Let $\sigma_{st}=(Z_{st},\mathcal{P}_{st})$ be the standard stability condition on the base curve. Suppose that $\sigma_{1}=(Z_{1},\mathcal{P}_{1})\in\Stab(p^{*}D^{b}(C)\otimes \mathcal{O}_{S}(-C_{0}))$ with $\mathcal{P}_{1}(0)=p^{*}\mathcal{P}_{st}(\phi_{1})\otimes \mathcal{O}_{S}(-C_{0})$ and $\sigma_{2}=(Z_{2},\mathcal{P}_{2})\in\Stab(p^{*}D^{b}(C))$ with $\mathcal{P}_{2}(0)=p^{*}\mathcal{P}_{st}(\phi_{2})$. Assume that $\sigma$ is a gluing pre-stability condition of $\sigma_{1}$ and $\sigma_{2}$, then \emph{gluing perversity of $\sigma$} is defined to be $\per(\sigma)=\phi_{1}-\phi_{2}$.
\end{def3.6}
\newtheorem*{prop3.7}{Proposition 3.7}
\begin{prop3.7}
Suppose $\sigma_{gl}$ is a gluing pre-stability condition. A $\widetilde{GL^{+}}(2,\mathbb{R})$-translate of $\sigma_{gl}$ has gluing perversity 1 if and only if $\per(\sigma_{gl})=1$
\end{prop3.7}
\begin{proof}
Suppose $\sigma_{gl}=(Z_{gl},\mathcal{P}_{gl})$ is a gluing pre-stability condition of $\sigma_{1}$ and $\sigma_{2}$, and $A=(M,f)\in\widetilde{GL^{+}}(2,\mathbb{R})$. If the heart of the bounded t-structure of $\sigma_{1}$ satisfies $\mathcal{P}_{1}(0)=p^{*}\mathcal{P}_{st}(\phi)\otimes\mathcal{O}(-C_{0})$ and the heart of the bounded t-structure of $\sigma_{2}$ satisfies $\mathcal{P}_{2}(0)=p^{*}\mathcal{P}_{st}(\psi)$, then $\per(\sigma_{gl}.A)=f^{-1}(\phi)-f^{-1}(\psi)$. $\per(\sigma_{gl})=\phi-\psi=1$ if and only if $\per(\sigma_{gl}.A)=f^{-1}(\phi)-f^{-1}(\psi)=1$ since $f$ is bijective.
\end{proof}
\newtheorem*{lem3.8}{Lemma 3.8}
\begin{lem3.8}
$\sigma_{1}$ and $\sigma_{2}$ satisfy the gluing property. Then $\per(\sigma)$ is not less than 1.
\end{lem3.8}
\begin{proof}
By Proposition 3.7, we can assume that $\sigma_{2}$ is the standard stability condition on $p^{*}D^{b}(S)$. Suppose that $\phi<1$ and $A_{1}=p^{*}\mathcal{P}_{st}(\phi,\phi+1]\otimes \mathcal{O}_{S}(-C_{0})$. It is enough to show that $\Hom(p^{*}\mathcal{P}_{st}(\phi,\phi+1])\otimes \mathcal{O}_{S}(-C_{0}),p^{*}\Coh C[i])\neq0$ for some $i\leq0$.\\
\ Recall that for all $q\in\frac{1}{\pi}\arctan\frac{1}{\mathbb{Z}}$ there is a line bundle $L$ such that $L\in \mathcal{P}_{st}(q)$. (For example, $L=\mathcal{O}_{C}(-n)$ with $q=\frac{1}{\pi}\arctan\frac{1}{n}$.) If we take $q\in(\phi-\lfloor\phi\rfloor,1)$, there is a line bundle $L\in \mathcal{P}_{st}(q)$ and we get $p^{*}L\otimes\mathcal{O}(-C_{0})[\lfloor\phi\rfloor]\in p^{*}\mathcal{P}_{st}(\phi,\phi+1]$. Hence, $\Hom(p^{*}L\otimes \mathcal{O}_{S}(-C_{0})[\lfloor\phi\rfloor],p^{*}L[\lfloor\phi\rfloor])\neq0$.
\end{proof}
\newtheorem*{thm3.9}{Theorem 3.9}
\begin{thm3.9}
On ruled surfaces, a gluing pre-stability condition $\sigma$ of $\widetilde{GL^{+}}(2,\mathbb{R})$-actions of the standard stability condition is a locally finite stability condition if and only if the gluing perversity of $\sigma$ is at least 1.
\end{thm3.9}
\begin{proof}
By Lemma 3.8, it would be sufficient to prove $\Hom(\mathcal{A}_{1},\mathcal{A}_{2}[i])$ for $i\leq0$ if $\phi=\per(\sigma)\geq 1$. By Proposition 3.7, we can assume $\mathcal{A}_{1}=p^{*}\mathcal{P}_{st}(\phi,\phi+1]\otimes\mathcal{O}_{S}(-C_{0})$ and $\mathcal{A}_{2}=p^{*}\mathcal{P}_{st}(0,1]$. Suppose that $F\in\mathcal{P}_{st}(\phi,\phi+1]$, $G\in\mathcal{P}_{st}(0,1]=\Coh C$ and $1\leq\phi$.\\
　$\Hom(p^{*}F\otimes \mathcal{O}_{S}(-C_{0}),p^{*}G[i])$\\
　$=\Hom(p^{*}F,p^{*}G\otimes \mathcal{O}_{S}(C_{0})[i])$\\
　$=\Hom(F,Rp_{*}(p^{*}G\otimes \mathcal{O}_{S}(C_{0})[i]))$\\
　$=\Hom(F,G\otimes Rp_{*}\mathcal{O}_{S}(C_{0})[i])$.\\
Since $Rp_{*}\mathcal{O}_{S}(C_{0})$ is a locally free sheaf, $G\otimes Rp_{*}\mathcal{O}_{S}(C_{0})[i]\in\mathcal{P}(i,i+1]$. Therefore, $\Hom(F,G\otimes Rp_{*}\mathcal{O}_{S}(C_{0})[i])=0$ by the phase of $F$ and $G\otimes Rp_{*}\mathcal{O}_{S}(C_{0})$. Then by Definition 3.2 (2), the image of $\sigma$ is discrete subgroup of $\mathbb{C}$. By \cite{CP} Proposition 3.5 (a), $\sigma$ is a Bridgeland stability condition.  Moreover, $\sigma$ is locally finite by \cite{Bri2} Lemma 4.4.
\end{proof}
In the above theorem, we declare all gluing stability conditions on ruled surfaces with base curve of positive genus. From now on, we mean a Bridgeland stability condition glued from $\widetilde{GL^{+}}(2,\mathbb{R})$-translates of stanard stability conditions on the base curve simply by a gluing stability conditions.
\newtheorem*{lem3.10}{Lemma 3.10}
\begin{lem3.10}
Suppose that $\sigma=(Z,\mathcal{A})$ is a gluing stability condition. Then\\
{\rm(1)} for any $f$, $\mathcal{O}_{f}$ and $\mathcal{O}_{f}(-C_{0})[1]$ are stable of the same phase in $\sigma$ respectively,\\
{\rm(2)} the phase of $\mathcal{O}_{f}$ is larger than the phase of $\mathcal{O}_{f}(-C_{0})[1]$,\\
{\rm(3)} if $\per(\sigma)=1$ skyscraper sheaves are strictly semistable of the same phase in $\sigma$, and also if $1<\per(\sigma)$ skyscraper sheaves are destabilised by $\mathcal{O}_{f}$ with $x\in f$. 
\end{lem3.10}
\begin{proof}
By Proposition 3.7, we can assume that $\sigma_{2}$ is the standard stability condition on $p^{*}D^{b}(S)$.\\
(1) Since $\mathcal{O}_{f}=p^{*}\mathcal{O}_{y}$ with $y=p(f)$, $\mathcal{O}_{f}$ is semistable of the same phase 1 for any $f$ by \cite{CP} Proposition 2.2 (3). Suppose that $\mathcal{F}$ is a subobject of $\mathcal{O}_{f}$ on $\mathcal{P}(1)$. $\mathcal{F}$ is also in $\mathcal{A}$. Hence, we have the following diagram in $\mathcal{P}(1)$.
$$
\begin{CD}
@. 0@. 0@. 0\\
@. @VVV @VVV @VVV @.\\
0@>>>\rho_{2}(\mathcal{F})@>>>\mathcal{F}@>>>\lambda_{1}(\mathcal{F})@>>>0\\
@. @VVV @VVV @VVV @.\\
0@>>>\rho_{2}(\mathcal{O}_{f})@>>>\mathcal{O}_{f}@>>>\lambda_{1}(\mathcal{O}_{f})@>>>0\\
\end{CD}
$$
Then $\mathcal{F}\simeq\rho_{2}(\mathcal{F})\subset \rho_{2}(\mathcal{O}_{f})=\mathcal{O}_{f}$ in $p^{*}D^{b}(C)$ by $\lambda_{1}(\mathcal{O}_{f})=0$. $\mathcal{O}_{f}$ is a minimal object in $p^{*}D^{b}(C)$. Hence, $\mathcal{F}$ is isomorphic to 0 or $\mathcal{O}_{f}$. $\mathcal{O}_{f}(-C_{0})[1]$ can be proved similarly.\\
(2) $\mathcal{O}_{f}=p^{*}\mathcal{O}_{y}$ with $y=p(f)$, $\mathcal{O}_{f}(-C_{0})[1]=p^{*}\mathcal{O}_{y}[1]\otimes \mathcal{O}_{S}(-C_{0})$ with $y=p(f)$. Since $\per(\sigma)\geq1$, the phase of $\mathcal{O}_{f}$ is larger than the phase of $\mathcal{O}_{f}(-C_{0})[1]$ by \cite{CP} Proposition 2.2 (3).\\
(3) If $\mathcal{O}_{x}$ is semistable of the phase $\phi$ we have the following in $\mathcal{A}[\lceil\phi\rceil-1]$. (c.f. \cite{CP} Lemma 2.1)\\
\centerline{$0\rightarrow\rho_{2}(\mathcal{O}_{x})\rightarrow \mathcal{O}_{x}\rightarrow\lambda_{1}(\mathcal{O}_{x})\rightarrow0$ exact.}\\
Since $\rho_{2}(\mathcal{O}_{x})=\mathcal{O}_{f}$ and $\lambda_{1}(\mathcal{O}_{x})=\mathcal{O}_{f}(-C_{0})[1]$, $\phi$ must be 1 by the phases, and hence if $1<\per(\sigma)$ $\mathcal{O}_{x}$ is destabilized by $\mathcal{O}_{f}$ with $x\in f$. Now we assume that $\per(\sigma)=1$. Since $\mathcal{O}_{f}\in \mathcal{P}(1)$ and $\mathcal{O}_{f}(-C_{0})\in \mathcal{P}(1)$, $\mathcal{O}_{x}$ is strictly semistable in $\sigma$ by extension closedness of $\mathcal{P}(1)$.
\end{proof}

\section{A destabilizing wall of skyscraper sheaves on ruled surfaces}
In this section, we describe a destabilizing wall of skyscraper sheaves on ruled surfaces. We start by the deformation theory of Bridgeland stability conditions.\\

For each $\sigma=(Z,\mathcal{P})\in\Stab_{\mathcal{N}}S$, define a function\\
\centerline{$||\cdot||_{\sigma}:\Hom(\mathcal{N}(S),\mathbb{C})\rightarrow[0,\infty)$}\\
by sending a group homomorphism $U:\mathcal{N}(S)\rightarrow\mathbb{C}$ to\\
\centerline{$||U||_{\sigma}=\sup\left\{\frac{|U(E)|}{|Z(E)|} \mid E\text{ semistable in }\sigma\right\}$}\\
Note that $||\cdot||_{\sigma}$ has all the properties of a norm on the complex vector space $\Hom(N(S),\mathbb{C})$. A norm of a finite dimensional vector space is unique up to equivalence. Hence, this norm is equivalent to the standard norm of the finite dimensional vector space  $\Hom(\mathcal{N}(S),\mathbb{C})$. If $\sigma=(Z,\mathcal{P})$ and $\tau=(W,\mathcal{Q})$ are stability conditions on a derived category $D^{b}(S)$ define\\
\centerline{$d(\mathcal{P},\mathcal{Q})=\sup\left\{|\phi^{+}_{\sigma}(E)-\phi^{+}_{\tau}(E)|,|\phi^{-}_{\sigma}(E)-\phi^{-}_{\tau}(E)| \mid 0\neq E\in D^{b}(S)\right\}$.}\\
It is a generalized metric on the space of slicings. Then an open basis of $\Stab_{\mathcal{N}}S$ consists of the following\\
\centerline{$B_{\mathcal{\epsilon}}(\sigma)=\left\{\tau=(W,\mathcal{Q})\in\Stab_{\mathcal{N}}S\mid||W-Z||_{\sigma}<\sin(\pi\epsilon) \text{ , } d(\mathcal{P},\mathcal{Q})<\epsilon\right\}.$}
\newtheorem*{prop4.1}{Proposition 4.1}
\begin{prop4.1}[\protect \cite{Bri1} Theorem 7.1]
Let $\sigma=(Z,\mathcal{P})$ be a numerical locally finite stability condition on a derived category $\mathcal{D}^{b}(S)$. Then there is an $\epsilon_{0}$ such that if $0<\epsilon<\epsilon_{0}$ and $W:\mathcal{N}(S)\rightarrow\mathcal{C}$ is a group homomorphism satisfying\\
\centerline{$|W(E)-Z(E)|<\sin(\pi\epsilon)|Z(E)|$}\\
for all $E\in\mathcal{D}^{b}(S)$ semistable in $\sigma$, then there is a locally finite stability condition $\tau=(W,\mathcal{Q})$ on $\mathcal{D}^{b}(S)$ with $d(\mathcal{P},\mathcal{Q})<\epsilon$.
\end{prop4.1}
The above $\mathcal{Q}$ is constructed as follows. A \emph{thin subcategory} of $\mathcal{D}^{b}(S)$ is a full subcategory of the form $\mathcal{P}((a,b))\subset\mathcal{D}^{b}(S)$ where $a$ and $b$ are real numbers with $0<b-a<1-2\epsilon$. Suppose $\psi(E)$ is the phase of $E$ on $W$. A nonzero object $E\in\mathcal{P}((a,b))$ is defined to be \emph{enveloped} by $\mathcal{P}((a,b))$ if $\mathcal{P}((a,b))$ is a thin subcategory satisfying $a+\epsilon\leq\psi(E)\leq b-\epsilon$. Then for each $\psi\in\mathbb{R}$ define $\mathcal{Q}(\psi)$ to be the full additive subcategory $\mathcal{D}^{b}(S)$ consisting of the zero objects of $\mathcal{D}^{b}(S)$ together with those object $E\in\mathcal{D}^{b}(S)$ which are $W$-semistable of phase $\psi$ in some thin enveloping subcategory $\mathcal{P}((a,b))$.\\

First, the following lemma plays an important role of the proof that gluing stability conditions with the gluing perversity 1 are a destabilizing wall of skyscraper sheaves.
\newtheorem*{lem4.2}{Lemma 4.2}
\begin{lem4.2}
Let $S$ be a ruled surface. Suppose that $\sigma_{gl}=(Z_{gl},P_{gl})$ is a gluing stability condition with the gluing perversity 1 on $S$. Then there is an $\epsilon_{0}>0$ such that if $0<\epsilon<\epsilon_{0}$ and $W:\mathcal{N}(S)\rightarrow\mathbb{C}$ is a group homomorphism satisfying\\
\centerline{the phase of $O_{f}(-C_{0})$ is greater than the phase of $O_{f}$, and}\\
\centerline{$|W(E)-Z(E)|<\sin(\pi\epsilon)|Z(F)|$}\\
for any $E\in D^{b}(S)$ semistable in $\sigma_{gl}$, then there is a unique locally finite Bridgeland stability condition $\tau=(W,\mathcal{Q})$ on $S$ with $d(\mathcal{P}_{gl},\mathcal{Q})<\epsilon$ satisfying that $O_{x}$ are stable of the same phase in $\tau$ for any $x\in S$.
\end{lem4.2}
\begin{proof}
By Proposition 3.7, we can assume that $\sigma_{2}$ is the standard stability condition on $p^{*}D^{b}(C)$. Then the phase of $\mathcal{O}_{x}$ is equal to 1. By the construction of $\mathcal{Q}, $we can construct the following slicing $\mathcal{Q}$ of $\tau$\\
\centerline{$\mathcal{Q}(\psi)=\left\{ F \middle| \begin{array}{l} \text{$F$ is enveloped by $\mathcal{P}_{gl}(a,b)$,}\\ \text{and semistable of phase $\psi$ in some $(W,\mathcal{P}_{gl}(a,b))$} \end{array} \right\}.$}\\
We show that $\mathcal{O}_{x}$ is a minimal object in $\mathcal{Q}(\psi)$. Since $\sigma_{gl}$ is discrete, we can take such an $\epsilon_{0}<\frac{1}{6}$ that\\
\centerline{$\mathcal{S}:=\{F\mid\re Z_{gl}(\mathcal{O}_{x})<\re Z_{gl}(F)<0$, $F\in \mathcal{P}_{gl}(1-2\epsilon,1+2\epsilon)\}\subset \mathcal{P}_{gl}(1)$.}\\
　　　　　\input{mtpic2.tex}\\
It is sufficient to show that $\mathcal{O}_{x}$ is stable in ($W$, $\mathcal{P}_{gl}(1-2\epsilon,1+2\epsilon)$). Suppose $\mathcal{O}_{x}$ is not stabile in $\mathcal{P}_{gl}(1-2\epsilon,1+2\epsilon)$. Then we can take $F$ a proper stable subobject of $\mathcal{O}_{x}$ in $\mathcal{P}_{gl}(1-2\epsilon,1+2\epsilon)$. We take an exact sequence in $\mathcal{P}_{gl}(1-2\epsilon,1+2\epsilon)$:\\
\centerline{$0\rightarrow F\rightarrow\mathcal{O}_{x}\rightarrow\mathcal{O}_{x}/F\rightarrow0$.}\\
We assume $F\notin\mathcal{P}_{gl}(1)$. Since $Z(\mathcal{O}_{x})=Z(F)+Z(\mathcal{O}_{x}/F)$, $\re Z(\mathcal{O}_{x})=\re Z(F)+\re Z(\mathcal{O}_{x}/F)$. Then we get $\re Z(\mathcal{O}_{x}/F)>0$ since $\re Z(F)\leq\re Z(\mathcal{O}_{x})\leq0$. This is contradictory to $\mathcal{O}_{x}/F\in\mathcal{P}_{gl}(1-2\epsilon,1+2\epsilon)$. Hence, we get $F\in\mathcal{P}_{gl}(1)$. We take $\alpha:F\hookrightarrow\mathcal{O}_{x}\rightarrow\mathcal{O}_{f}(-C_{0})[1]$.
\begin{itemize}
\item If $\alpha=0$, there exists a morphism $F\rightarrow\mathcal{O}_{f}$. 
$$
\begin{CD}
F@>>>\mathcal{O}_{x}@>>>\mathcal{O}_{x}/F@>>>F[1]\\
@VVV @VVV @. @.\\
\mathcal{O}_{f}@>>>\mathcal{O}_{x}@>>>\mathcal{O}_{f}(-C_{0})[1]@>>>\mathcal{O}_{f}[1]\\
\end{CD}
$$
Since $\mathcal{O}_{f}$ is a minimal object in $\mathcal{P}_{gl}(1)$, we get $F\simeq\mathcal{O}_{f}$.
\item If $\alpha\neq0$, $\alpha$ is surjective. Moreover, we get $\ker\alpha\simeq 0$ since $F$ is stable in $\mathcal{P}_{gl}(1)$. Hence $\alpha$ is ismorphism. So $F\simeq\mathcal{O}_{f}(-C_{0})[1]$. Since $\Hom(F,\mathcal{O}_{x})=\Hom(\mathcal{O}_{f}(-C_{0})[1],\mathcal{O}_{x})=0$, then this is contradictory to $F\subset\mathcal{O}_{x}$.
\end{itemize}
Hence, we get $F\simeq\mathcal{O}_{f}$. Since $W(\mathcal{O}_{x})=W(\mathcal{O}_{f})+W(\mathcal{O}_{f}(-C_{0})[1])$ and $\psi(\mathcal{O}_{f})<\psi(\mathcal{O}_{f}(-C_{0})[1])$, $\psi(\mathcal{O}_{f})<\psi(\mathcal{O}_{x})=\psi$. Namely, $\mathcal{O}_{x}$ is stable in ($W$,$\mathcal{P}_{gl}(1-2\epsilon,1+2\epsilon)$).
\end{proof}
Second, the set of gluing stability conditions are connected submanifold of $\Stab_{\mathcal{N}}S$. We prove the following lemma.
\newtheorem*{lem4.3}{Lemma 4.3}
\begin{lem4.3}
Let $S_{gl,p}$ be the set of gluing stability conditions with gluing perversity p. $S_{gl,1}$ is connected submanifold of $\Stab_{\mathcal{N}}S$ with real dimension 7. Moreover, $S_{gl}:=\bigcup_{p}S_{gl,p}$ is also a submanifold with real dimension 8, especially the subset of full components.
\end{lem4.3}
\begin{proof}
We show that the action of $\widetilde{GL^{+}}(2,\mathbb{R})$ on $S_{gl,1}$ is free. Suppose $\sigma_{gl}\in S_{gl,1,st}$ and $A=(M,f)\in\widetilde{GL^{+}}(2,\mathbb{R})$. If $\sigma_{gl}.A=\sigma_{gl}$, then we get\\
\centerline{$M^{-1}(Z_{gl}(\mathcal{O}_{S}))=Z_{gl}(\mathcal{O}_{S})$}\\
and\\
\centerline{$M^{-1}(Z_{gl}(\mathcal{O}_{f}))=Z_{gl}(\mathcal{O}_{f})$.}\\
By Proposition 3.5, $Z_{gl}(\mathcal{O}_{S})=i$ and $Z_{gl}(\mathcal{O}_{f})=-1$. Hence, $M$ is the identity matrix by comparison of both values of central charges. $f=\id$ can be get by the comparison of both hearts of the bounded t-structures. Suppose that $S_{gl,1,st}$ consists of the element of $S_{gl,1}$ that $\sigma_{2}$ is the standard stability condition on $p^{*}D^{b}(C)$. Then by \cite{Bri1} Theorem 9.1, $S_{gl,1,st}\simeq\{(\begin{pmatrix}a&&b\\0&&d\end{pmatrix},f)\mid a>0,b\in\mathbb{R},d>0\text{ and }f(0)=0\}$. Especially, $S_{gl,1,st}$ is a connected submanifold with real dimension 3 since $\pr_{1}$ is a local homeomorphism. Hence, $S_{gl,p}$ is connected submanifold of $\Stab_{\mathcal{N}}S$ with real dimension 7. We can prove in the case of $S_{gl}$ similarly.
\end{proof}
Finally, we describe a concrete description between geometric stability conditions and gluing stability conditions on the stability space. This is the end of the proof of Theorem 1.4.
\newtheorem*{thm4.4}{Theorem 4.4}
\begin{thm4.4}
Let $S_{geom}$ be the set of geometric stability conditions on $S$. Suppose that $A=(\begin{pmatrix}a&&\frac{1}{2}a\deg\mathcal{E}\\0&&a\end{pmatrix}^{-1},f)\in\widetilde{GL^{+}}(2,\mathbb{R})$ with $a<0$. Then $\partial\overline{S_{geom}}\cap S_{gl,1}$ is the set of $\widetilde{GL^{+}}(2,\mathbb{R})$-translates of a stability condition glued from $\sigma_{st}.A$ and $\sigma_{st}$. 
\end{thm4.4}
\begin{proof}
We can assume that $\sigma_{gl}=(Z_{gl},P_{gl})$ is a gluing stability condition that $\sigma_{2}$ is a standard stabilty condition. It is sufficient to show that $Z_{gl}=M^{-1}\exp(B+i\omega)$ if and only if $Z_{gl}=\begin{pmatrix}a&&\frac{1}{2}a\deg\mathcal{E}\\0&&a\end{pmatrix}Z_{st}\circ\lambda_{1}+Z_{st}\circ\rho_{2}$ with $a<0$.\\
~　Let $M^{-1}=\begin{pmatrix}\alpha&&\beta\\\gamma&&\delta\end{pmatrix}$, $B=xC_{0}+yf$ and $\omega=zC_{0}+wf$. We denote $I=\frac{1}{2}\alpha\{(x^{2}-z^{2})\deg\mathcal{E}+2(xy-zw)\}+\beta\{xz\deg\mathcal{E}+(yz+xw)\}$ and $J=\frac{1}{2}\gamma\{(x^{2}-z^{2})\deg\mathcal{E}+2(xy-zw)\}+\delta\{xz\deg\mathcal{E}+(yz+xw)\}$ Then\\
　$\exp(B+i\omega)$\\
　$=(1,x+iz,y+iw,\frac{1}{2}\{(x^{2}-z^{2})\deg\mathcal{E}+2(xy-zw)\}+i\{xz\deg\mathcal{E}+(yz+xw)\})$,\\
　$M^{-1}\exp(B+i\omega)$\\
　$=(\alpha+i\gamma,\{(\alpha x+\beta z)+i(\gamma x+\delta z)\}C_{0}+\{(\alpha y+\beta w)+i(\gamma y+\delta w)\}f, I+iJ)$.\\
We compare it to Proposition 3.5. Recall that $\sigma_{gl}$ has gluing perversity 1. So $a<0$ and $c=0$. Then\\
　$\pr_{1}(Z_{gl})=(1-a, -C_{0}+[\{\frac{1}{2}\deg\mathcal{E}(a+1)-b\}+i(1-d)]f,-i)$\\
From $\alpha+i\gamma=1-a$, we get $\alpha=1-a$ and $\gamma=0$. Then we get $z=0$ from $\gamma x+\delta z=0$ since $\det M=\alpha\delta\neq0$. And then we get $x=\frac{1}{a-1}$ from $\alpha x+\beta z=-1$. And then we get $a=d$ from $J=\delta xw=-1$ and $\gamma y+\delta w=1-d$. From $I=-\frac{1}{2}(\frac{1}{a-1}\deg\mathcal{E}+2y)+\beta\frac{1}{a-1}w=0$ and $\alpha y+\beta w=(1-a)y+\beta w=\frac{1}{2}\deg\mathcal{E}(a+1)-b$, we get $b=\frac{1}{2}a\deg\mathcal{E}$.
\end{proof}
The set of gluing stability conditions is a submanifold of the full stability space (Lemma 4.3). Lemma 3.10 (3) and Lemma 4.2 suggest that the set of gluing stability conditions neighbors on the set of stability conditions such that skyscraper sheaves are stable of the same phase on the stability space. Especially, the set of gluing stability conditions with the gluing perversity 1 is a part of destabilizing wall of skyscraper sheaves. In addition, the boundary of the set of geometric stability conditions only contacts the destabilizing wall (Theorem 4.4). The following picture of $\Stab_{\mathcal{N}}S$ is convenient for understanding.\\
\\
\input{mtpic1.tex}\\
\newtheorem*{rem45}{Remark 4.5}
\begin{rem45}
Let $\overline{\mathcal{M}^{\sigma}}([\mathcal{O}_{x}])$ be the variety of S-equivalent classes of objects $E\in\mathcal{P}(\phi(\mathcal{O}_{x}))$. 
\begin{itemize}
\item If $\sigma$ is a geometric stability condition, then $\overline{\mathcal{M}^{\sigma}}([\mathcal{O}_{x}])\simeq S$.
\item If $\sigma$ is a gluing stability condition with gluing perversity $1$, then $\overline{\mathcal{M}^{\sigma}}([\mathcal{O}_{x}])\simeq C$.
\item If $\sigma$ is a gluing stability condition with gluing perversity $>1$, then $\overline{\mathcal{M}^{\sigma}}([\mathcal{O}_{x}])$ is empty.
\end{itemize}
\end{rem45}

\end{document}

%% file: mtpic2.tex
{\unitlength 0.1in%
\begin{picture}(14.3000,13.3000)(8.0000,-20.0000)%
%
\special{pn 8}%
\special{pa 800 1400}%
\special{pa 2200 1400}%
\special{fp}%
\special{sh 1}%
\special{pa 2200 1400}%
\special{pa 2133 1380}%
\special{pa 2147 1400}%
\special{pa 2133 1420}%
\special{pa 2200 1400}%
\special{fp}%
\special{pa 2000 2000}%
\special{pa 2000 800}%
\special{fp}%
\special{sh 1}%
\special{pa 2000 800}%
\special{pa 1980 867}%
\special{pa 2000 853}%
\special{pa 2020 867}%
\special{pa 2000 800}%
\special{fp}%
\put(19.8000,-8.0000){\makebox(0,0)[lb]{Im}}%
\put(22.3000,-14.4000){\makebox(0,0)[lb]{Re}}%
%
\special{pn 8}%
\special{pa 2000 1400}%
\special{pa 800 1000}%
\special{fp}%
\special{pa 820 1800}%
\special{pa 2000 1400}%
\special{fp}%
\special{pa 1000 1070}%
\special{pa 1000 1750}%
\special{fp}%
%
\special{pn 4}%
\special{pa 1360 1400}%
\special{pa 1030 1730}%
\special{fp}%
\special{pa 1300 1400}%
\special{pa 1000 1700}%
\special{fp}%
\special{pa 1240 1400}%
\special{pa 1000 1640}%
\special{fp}%
\special{pa 1180 1400}%
\special{pa 1000 1580}%
\special{fp}%
\special{pa 1120 1400}%
\special{pa 1000 1520}%
\special{fp}%
\special{pa 1060 1400}%
\special{pa 1000 1460}%
\special{fp}%
\special{pa 1420 1400}%
\special{pa 1130 1690}%
\special{fp}%
\special{pa 1480 1400}%
\special{pa 1220 1660}%
\special{fp}%
\special{pa 1540 1400}%
\special{pa 1310 1630}%
\special{fp}%
\special{pa 1600 1400}%
\special{pa 1400 1600}%
\special{fp}%
\special{pa 1660 1400}%
\special{pa 1490 1570}%
\special{fp}%
\special{pa 1720 1400}%
\special{pa 1580 1540}%
\special{fp}%
\special{pa 1780 1400}%
\special{pa 1670 1510}%
\special{fp}%
\special{pa 1840 1400}%
\special{pa 1760 1480}%
\special{fp}%
\special{pa 1900 1400}%
\special{pa 1850 1450}%
\special{fp}%
%
\special{pn 4}%
\special{pa 1290 1170}%
\special{pa 1060 1400}%
\special{fp}%
\special{pa 1250 1150}%
\special{pa 1010 1390}%
\special{fp}%
\special{pa 1200 1140}%
\special{pa 1000 1340}%
\special{fp}%
\special{pa 1160 1120}%
\special{pa 1000 1280}%
\special{fp}%
\special{pa 1110 1110}%
\special{pa 1000 1220}%
\special{fp}%
\special{pa 1070 1090}%
\special{pa 1000 1160}%
\special{fp}%
\special{pa 1020 1080}%
\special{pa 1000 1100}%
\special{fp}%
\special{pa 1340 1180}%
\special{pa 1120 1400}%
\special{fp}%
\special{pa 1380 1200}%
\special{pa 1180 1400}%
\special{fp}%
\special{pa 1430 1210}%
\special{pa 1240 1400}%
\special{fp}%
\special{pa 1470 1230}%
\special{pa 1300 1400}%
\special{fp}%
\special{pa 1520 1240}%
\special{pa 1360 1400}%
\special{fp}%
\special{pa 1560 1260}%
\special{pa 1420 1400}%
\special{fp}%
\special{pa 1610 1270}%
\special{pa 1480 1400}%
\special{fp}%
\special{pa 1650 1290}%
\special{pa 1540 1400}%
\special{fp}%
\special{pa 1700 1300}%
\special{pa 1600 1400}%
\special{fp}%
\special{pa 1740 1320}%
\special{pa 1660 1400}%
\special{fp}%
\special{pa 1790 1330}%
\special{pa 1720 1400}%
\special{fp}%
\special{pa 1830 1350}%
\special{pa 1780 1400}%
\special{fp}%
\special{pa 1880 1360}%
\special{pa 1840 1400}%
\special{fp}%
\put(20.2000,-15.4000){\makebox(0,0)[lb]{O}}%
\put(15.0000,-11.9000){\makebox(0,0)[lb]{$\mathcal{S}$: There are objects only on real axis.}}%
\end{picture}}%

%% file: mtpic1.tex
{\unitlength 0.1in%
\begin{picture}( 27.8000, 14.0000)( 12.9000,-18.0000)%
%
\special{pn 8}%
\special{ar 2800 1200 1200 600  1.5946014  1.5707963}%
%
\special{pn 8}%
\special{pa 2610 610}%
\special{pa 3010 1790}%
\special{fp}%
%
\special{pn 8}%
\special{ar 3410 1210 562 562  1.7761917  4.5679474}%
%
\special{pn 4}%
\special{pa 2990 610}%
\special{pa 2710 890}%
\special{fp}%
\special{pa 3100 620}%
\special{pa 2740 980}%
\special{fp}%
\special{pa 2950 890}%
\special{pa 2770 1070}%
\special{fp}%
\special{pa 2860 1100}%
\special{pa 2800 1160}%
\special{fp}%
\special{pa 2850 1230}%
\special{pa 2830 1250}%
\special{fp}%
\special{pa 3200 640}%
\special{pa 3090 750}%
\special{fp}%
\special{pa 2880 600}%
\special{pa 2680 800}%
\special{fp}%
\special{pa 2760 600}%
\special{pa 2650 710}%
\special{fp}%
%
\special{pn 4}%
\special{pa 3110 1690}%
\special{pa 3020 1780}%
\special{fp}%
\special{pa 3190 1730}%
\special{pa 3150 1770}%
\special{fp}%
\special{pa 3040 1640}%
\special{pa 2980 1700}%
\special{fp}%
\special{pa 2980 1580}%
\special{pa 2950 1610}%
\special{fp}%
%
\special{pn 4}%
\special{pa 3690 810}%
\special{pa 2960 1540}%
\special{fp}%
\special{pa 3730 830}%
\special{pa 2990 1570}%
\special{fp}%
\special{pa 3770 850}%
\special{pa 3020 1600}%
\special{fp}%
\special{pa 3800 880}%
\special{pa 3050 1630}%
\special{fp}%
\special{pa 3830 910}%
\special{pa 3080 1660}%
\special{fp}%
\special{pa 3860 940}%
\special{pa 3120 1680}%
\special{fp}%
\special{pa 3890 970}%
\special{pa 3160 1700}%
\special{fp}%
\special{pa 3920 1000}%
\special{pa 3200 1720}%
\special{fp}%
\special{pa 3940 1040}%
\special{pa 3240 1740}%
\special{fp}%
\special{pa 3960 1080}%
\special{pa 3300 1740}%
\special{fp}%
\special{pa 3980 1120}%
\special{pa 3380 1720}%
\special{fp}%
\special{pa 3990 1170}%
\special{pa 3460 1700}%
\special{fp}%
\special{pa 3990 1230}%
\special{pa 3570 1650}%
\special{fp}%
\special{pa 3970 1310}%
\special{pa 3690 1590}%
\special{fp}%
\special{pa 3650 790}%
\special{pa 2940 1500}%
\special{fp}%
\special{pa 3610 770}%
\special{pa 2920 1460}%
\special{fp}%
\special{pa 3570 750}%
\special{pa 2900 1420}%
\special{fp}%
\special{pa 3530 730}%
\special{pa 2880 1380}%
\special{fp}%
\special{pa 3490 710}%
\special{pa 2870 1330}%
\special{fp}%
\special{pa 3440 700}%
\special{pa 2860 1280}%
\special{fp}%
\special{pa 3400 680}%
\special{pa 2850 1230}%
\special{fp}%
\special{pa 3350 670}%
\special{pa 2850 1170}%
\special{fp}%
\special{pa 3300 660}%
\special{pa 2860 1100}%
\special{fp}%
\special{pa 3210 690}%
\special{pa 2890 1010}%
\special{fp}%
%
\special{pn 8}%
\special{pa 3250 490}%
\special{pa 2870 760}%
\special{fp}%
\special{sh 1}%
\special{pa 2870 760}%
\special{pa 2936 738}%
\special{pa 2913 729}%
\special{pa 2913 705}%
\special{pa 2870 760}%
\special{fp}%
\put(33.0000,-5.3000){\makebox(0,0)[lb]{All point sheaves are stable of the same phase.}}%
%
\special{pn 8}%
\special{pa 4010 1600}%
\special{pa 3620 1400}%
\special{fp}%
\special{sh 1}%
\special{pa 3620 1400}%
\special{pa 3670 1448}%
\special{pa 3667 1424}%
\special{pa 3688 1413}%
\special{pa 3620 1400}%
\special{fp}%
\put(40.7000,-16.6000){\makebox(0,0)[lb]{geometric stability conditions}}%
\put(12.9000,-15.0000){\makebox(0,0)[lb]{gluing stability conditions}}%
\put(30.7000,-19.3000){\makebox(0,0)[lb]{gluing perversity 1}}%
\end{picture}}%